\documentclass[12pt]{article}
\usepackage[final]{epsfig}
\usepackage{graphics}
\usepackage{amsmath}
\usepackage{amsfonts}
\usepackage{latexsym}
\usepackage{amssymb}
\usepackage{graphicx}
\usepackage{epstopdf}
\newtheorem{lemma}{Lemma}[section]

\newtheorem{remark}[lemma]{Remark}
\newtheorem{example}[lemma]{Example}
\newtheorem{theorem}{Theorem}

\begin{document}
\newcommand{\eps}{{\varepsilon}}
\newcommand{\proofend}{$\Box$\bigskip}
\newcommand{\C}{{\mathbf C}}
\newcommand{\Q}{{\mathbf Q}}
\newcommand{\R}{{\mathbf R}}
\newcommand{\Z}{{\mathbf Z}}
\newcommand{\RP}{{\mathbf {RP}}}
\newcommand{\CP}{{\mathbf {CP}}}

\title {Converse Sturm-Hurwitz-Kellogg theorem and related results}
\author{Serge Tabachnikov\thanks{
Department of Mathematics,
Pennsylvania State University, University Park, PA 16802, USA;
e-mail: \tt{tabachni@math.psu.edu}}
}
\date{\today}
\maketitle
\begin{abstract} We prove that  if $V^n$ is a Chebyshev system on the circle and $f(x)$ is a continuous function with at least $n+1$ sign changes then there exists an orientation preserving diffeomorphism of $S^1$ that takes $f$ to a function $L^2$-orthogonal to $V$. We also prove that if $f(x)$ is a function on the real projective line with at least four sign changes then there exists an orientation preserving diffeomorphism of $\RP^1$ that takes $f$ to the Schwarzian derivative of a function on $\RP^1$. We show that the space of piece-wise constant functions on an interval with values $\pm1$ and at most $n+1$ intervals of constant sign is homeomorphic to $n$-dimensional sphere.

\end{abstract}

\vskip 10 mm

{\hfill To V. I. Arnold for his 70th birthday}

\bigskip

\section{Introduction and formulation of results} \label{intro}

The classic four vertex theorem asserts that the curvature of a plane oval (strictly convex smooth closed curve) has at least four extrema. Discovered about 100 years ago by S. Mukhopadhyaya, this theorem and its numerous generalizations and refinements continue to attract attention up to this day; see \cite{O-T} for a sampler. 

One such result is the converse four vertex theorem proved by  Gluck for strictly convex, and by Dahlberg for general curves  \cite{G1,Da}: a periodic function having at least two local minima and two local maxima is the curvature function of a simple closed plane curve. See \cite{G2} for a very well written survey.

The radius of curvature $\rho(\alpha)$ of an oval, considered as a function of the direction of the tangent line to the curve, is $L^2$-orthogonal to the first harmonics:
$$
\int_0^{2\pi} \rho(\alpha) \cos\alpha\ d\alpha=\int_0^{2\pi} \rho(\alpha) \sin\alpha\ d\alpha=0.
$$
Such a function must have at least four critical points. The converse four vertex theorem can be restated as follows: {\it if a function $\rho(\alpha)$ has at least two local minima and two local maxima then there is a diffeomorphism $\varphi$ of the circle such that the function $\rho(\varphi(\alpha))$ is $L^2$-orthogonal to the first harmonics.}

Our first result is the following generalization. 

A Chebyshev system is an $n$-dimensional\footnote{Where $n$ is odd. One can define a Chebyshev system on a segment as well, and then there is no restriction on the parity of its dimension.} space $V$ of functions on the circle $S^1=\R/2\pi \Z$ such that every non-zero function from $V$ has at most $n-1$ zeros (counted with multiplicities). According to the Sturm-Hurwitz-Kellogg theorem, if a smooth function\footnote{Smoothness is not needed; one can work with finitely differentiable or continuous functions.}  on $S^1$ is $L^2$-orthogonal to a Chebyshev system $V^n$  then this function has at least $n+1$ sign changes; see, e.g., \cite{O-T}. In particular, a function orthogonal to $\{1,\cos\alpha,\sin\alpha\}$ has at lest four zeros; applied to the derivative of the radius of curvature of an oval, this implies the four vertex theorem.

We prove the next converse Sturm-Hurwitz-Kellogg theorem.

\begin{theorem} \label{main}
Let $V^n$ be a Chebyshev system on $S^1$. If $f(x)$ is a continuous function on $S^1$ with at least $n+1$ sign changes then there exists an orientation preserving diffeomorphism  $\varphi: S^1\to S^1$ such that $f(\varphi(x))$ is $L^2$-orthogonal to $V$.
\end{theorem}

Our strategy of proof is that of Gluck  \cite{G1,G2} which we illustrate by the following simplest case of the above theorem.

\begin{example} \label{exa}
{\rm Let $f(x)$ be a continuous  function on $S^1$ that has both positive and negative values. One claims that there exists an orientation preserving diffeomorphism  $\varphi: S^1\to S^1$ such that $f(\varphi(x))$ has zero average value:
$$
\int_0^{2\pi} f(\varphi(x))\ dx =0.
$$ 
Of course, this is obvious, but we shall describe an argument that exemplifies the method 
of proof of Theorem \ref{main} and other results of this paper.

\paragraph{\it Step 1.} Let $h(x)$ be the step function that takes value $1$ on $[0,\pi)$ and $-1$ on $[\pi,2\pi)$. This step function has zero average value.

\paragraph{\it Step 2.} Since $f(x)$ changes sign, there is a number $c\neq 0$ such that $f$ assumes both values $\pm c$. Scaling $f$, assume that $c=1$ and that $f(x_1)=1, f(x_2)=-1$.   For every $\eps >0$, there exists a diffeomorphism $\varphi \in {\rm Diff}_+ (S^1)$ which stretches neighborhoods of the points  $x_1$ and $x_2$ so that $\varphi^*(f)$ is $\eps$-close in measure to $h$.

\paragraph{\it Step 3.} For a sufficiently small real $\alpha$, consider an orientation preserving diffeomorphism $\psi_{\alpha} \in {\rm Diff}_+ (S^1)$ that fixes $0$ and stretches the interval $[0,\pi]$ to $[0,\pi+\alpha]$. We assume that  the dependence of  $\psi_{\alpha}$ on $\alpha$ is smooth. The correspondence $\alpha\mapsto \psi_{\alpha}$ is a map of an interval $I$ to the group ${\rm Diff}_+ (S^1)$. Consider the function 
$$
F(\alpha) =\int_0^{2\pi} (\psi_{\alpha}^*)(h)(x)\ dx.
$$
One has: $F(0)=0$ and $F'(0)\neq 0$. In particular, making the interval $I$ smaller, if needed, $F$ has opposite signs at the end points of $I$.

\paragraph{\it Step 4.} Finally, replace $h$ in the definition of $F$ by the function $\varphi^*(f)$ from Step 2. If $\eps$ is small enough, the resulting function $\bar F: I \to \R$ still has opposite signs at the end points of $I$, hence there exists $\alpha$ such that $\bar F(\alpha)=0$. Thus the function  $\psi_{\alpha}^* (\varphi^*(f))$ has zero average.
}
\end{example}

\begin{remark} \label{differential}
{\rm An object invariantly related to a function is its differential $df=f'(x) dx$ (rather than the derivative). If $\lambda$ is a differential 1-form on $S^1$ and
$$
\int_0^{2\pi} \lambda =0
$$
then $\lambda$ has sign changes, but the converse does not hold since 
$$
\int_0^{2\pi} \varphi^*(\lambda)= \int_0^{2\pi} \lambda 
$$
for every $\varphi\in {\rm Diff}_+ (S^1)$. This explains why we deal with a function, rather than a differential 1-form.
}
\end{remark}

Another, rather recent, four vertex-type theorem is due to E. Ghys: the Schwarzian derivative of  a diffeomorphism of the real projective line has at least four zeros. Choose an affine coordinate $x$ on $\RP^1$ and let $f(x)$ be a diffeomorphism. Then the Schwarzian derivative $S(f)$ is given by the formula
$$
S(f)=\frac{f'''}{f'} - \frac{3}{2} \left( \frac{f''}{f'} \right)^2;
$$
it measures the failure of $f$ to preserve the projective structure; see \cite{O-T}.

We prove a converse theorem.

\begin{theorem} \label{Schw}
If $f(x)$ is a smooth function on $\RP^1$ with at least four sign changes then there exists an orientation preserving diffeomorphisms of the projective line  $\varphi$ and $g(x)$ such that $\varphi^*(f)=S(g)$.
\end{theorem}

\begin{remark} \label{remSc}
{\rm The invariant meaning of the Schwarzian is not a function but rather a quadratic differential, see, e.g., \cite{O-T} for a detailed discussion:
$$
S(f)=\left(\frac{f'''}{f'} - \frac{3}{2} \left( \frac{f''}{f'} \right)^2\right) dx^2.
$$
Similarly to Remark \ref{differential}, the property of a quadratic differential on $\RP^1$ to be the Schwarzian derivative of a diffeomorphism is invariant under the action of the group ${\rm Diff}(\RP^1)$. 
}
\end{remark}

\section{Proof of the converse Sturm-Hurwitz-Kellogg theorem} \label{SHK}

The proof consists of the same four steps as in Example \ref{exa}. 

\paragraph{\it Step 1.}
\begin{lemma} \label{orth}
There exists a piece-wise constant function on $S^1$ with values $\pm 1$ and exactly $n+1$ intervals of constant sign which is $L^2$-orthogonal to $V$.
\end{lemma}

\paragraph{\bf Proof} (suggested by D. Khavinson). 
Extend $V^n$ to a larger Chebyshev system $W^{n+2}$ and pick $f\in W-V$. Consider $g$, the best $L^1$ approximation of $f$ by a function in $V$. The function $g$ exists since $V$ is finite dimensional. 

Since $W$ is a Chebyshev system, $f-g$ has at most $n+1$ intervals of constant sign (obviously, $f-g\neq 0$). Let $I_k$ be these intervals, and let $h$ be the function that has alternating values $\pm 1$ on the intervals $I_k$. Since $g$ is best approximation of $f$, one has the Lagrange multipliers condition:
\begin{equation} \label{Lagr}
\left.\frac{d}{d \eps}\right|_{\eps = 0}  \left( \int_0^{2\pi} |(f-g)(x)+\eps v(x)|\ dx\right)  = 0
\end{equation}
for every $v(x) \in V$. It follows from (\ref{Lagr}) that 
$$
0=\sum_k (-1)^k \int_{I_k} v(x)\ dx =  \int_0^{2\pi} h(x) v(x)\ dx,
$$
that is, $h$ is orthogonal to $V$.

By the Sturm-Hurwitz-Kellog theorem, $ h$ has at least $n+1$ sign changes (Proof, for completeness: if not, one can find a function from $V$ with the same intervals of constant sign as $h$; such a function cannot be orthogonal to $h$).
\proofend

\paragraph{\it Step 2.} Since $f(x)$ changes sign at least $n+1$  times, there is a non-zero constant $c$ such that $f$ takes the alternating values $\pm c$ at points, say, $x_0,\dots,x_n$. Multiplying $f$ by a constant, assume that $c=1$. 

Let $h(x)$ be the function from Lemma \ref{orth}. For every $\eps >0$, there exists a diffeomorphism $\varphi \in {\rm Diff}_+ (S^1)$ which stretches neighborhoods of the points  $x_0,\dots,x_n$ so that the function $\varphi^*(f)$ is $\eps$-close in measure to $h$.

\paragraph{\it Step 3.} Consider the function $h(x)$ and let $[0,x_1], [x_1,x_2],\dots,[x_n,2\pi]$ be its intervals of constant sign. For $\alpha=(\alpha_1,...,\alpha_n)$, consider an orientation preserving diffeomorphism $\psi_{\alpha} \in {\rm Diff}_+ (S^1)$ that stretches the intervals $[x_i,x_{i+1}]$ so that  point $x_i$ goes to $x_i + \alpha_i$ and which fixes $0$. We assume that each $|\alpha_i|$ is sufficiently small and that the dependence of  $\psi_{\alpha}$ on $\alpha$ is smooth. The correspondence $\alpha\mapsto \psi_{\alpha}$ is a map of an $n$-dimensional disc $D^n$ to ${\rm Diff}_+ (S^1)$.

The formula $F(\alpha) (g)=\langle \psi_{\alpha}^*(h),g\rangle$ defines a smooth map 
$D\to V^*$ that takes the origin to the origin (the scalar product is understood in the $L^2$ sense). 

\begin{lemma} \label{surj}
The differential $dF$ is non-degenerate at the origin.
\end{lemma}

\paragraph{\bf Proof.} Let $g_1,\dots,g_n$ be a basis of $V$.
 We want to prove that the matrix 
$$
c_{ij}=\left.\frac{\partial F(\alpha) (g_i)}{\partial \alpha_j}\right|_{\alpha=0},\quad i,j=1,\dots,n
$$
is non-degenerate.  One has: 
$$
F(\alpha) (g) = \sum_{k=0}^{n} (-1)^k \int_{x_k+\alpha_k}^{x_{k+1}+\alpha_{k+1}} g(x)\ dx
$$
where we assume that $x_0=0, x_{n+1}=2\pi, \alpha_0=\alpha_{n+1}=0$. It follows that 
$c_{ij}=2(-1)^{j+1} g_i (x_j)$, and it suffices to show that the matrix $g_i (x_j)$ is non-degenerate. This is indeed a fundamental property of Chebyshev systems, see \cite{K-N} (Proof, for completeness: if $c=(c_1,\dots,c_n)$ is a non-zero vector such that $\sum c_i g_i (x_j)=0$ for each $j$ then the function $\sum c_i g_i(x)$ has $n$ zeros, which contradicts the definition of Chebyshev systems).
\proofend

\paragraph{\it Step 4.} It follows from Lemma \ref{surj} that there exists $\delta>0$ such that the map $F$, restricted to the cube $D^n$ given by the conditions $|\alpha_i|<\delta,\ i=1,\dots,n$, has degree one, and the hypersurface $F(\partial D)$ has the rotation number one with respect to the origin in $V^*$.

Now replace $h$ in the definition of the map $F$ by the function $\varphi^*(f)$ from Step 2, and denote the new map by $\bar F: D^n \to V^*$. We shall be done if we  show that there exists $\alpha$ such that $\bar F(\alpha)=0$. Indeed, if $\eps$ is small enough then $\bar F(\partial D)$ still has rotation number one with respect to the origin in $V^*$, and therefore $\bar F(D)$ contains the origin. 
\proofend

\section{Digression: the space of step functions with values $\pm 1$ on an interval} \label{step}
An extension of Lemma \ref{orth}  to the case when  $V$ is not assumed to be a Chebyshev system  is the following Hobby--Rice theorem \cite{H-R}, see also   \cite{Pi,Tot}.

\begin{theorem} \label{alt}
Let $V$ be an $n$-dimensional subspace in $L^1([0,1])$. Then there exists a piece-wise constant function on $I$ with values $\pm 1$ and at most $n+1$ intervals of constant sign which is $L^2$-orthogonal to $V$.
\end{theorem}

\paragraph{\bf Proof} (\cite{H-R,Pi}).
Let $x=(x_0,x_1,\dots,x_n),\ \sum_i x_i^2=1$, be a point of the sphere $S^n$. Assign to $x$ the partition of $[0,1]$ on the intervals of consecutive lengths $x_0^2,\dots,x_n^2$ and the piece-wise constant function $h_x$ with value equal to  sign $x_i$ on the respective interval. We obtain a map  $F:S_n\to V^*$  given by the formula: 
$$
\langle F(x), g\rangle= \int_0^1 h_x(t) g(t)\ dt.
$$
This map is odd: $F(-x)=-F(x)$, and it follows from the Borsuk-Ulam theorem (see e.g., \cite{Ma}) that $F(x)=0$ for some $x\in S_n$. Thus $h_x$ is orthogonal to $V$.
\proofend

From the point of view of topology, it is interesting to consider the space $S_n \subset L^1([0,1])$  of piece-wise constant function on $[0,1]$ with values $\pm 1$ and at most $n+1$ intervals of constant sign. We complement the proof of Theorem \ref{alt} with the following result.

\begin{theorem} \label{sphere}
$S_n$ is homeomorphic to $n$-dimensional sphere.
\end{theorem}

\paragraph{\bf Proof.} We give $S_n$ the structure of a finite cell complex with two cells in every dimension $0,1,\dots,n$ and prove, by induction on $n$, that $S_n$ is homeomorphic to $S^n$. For $n=0$, the set $S_0$ consists of two constant functions with values $+1$ or $-1$ and is homeomorphic to $S^0$.

Let $\Delta^n = \{x=(x_0,\dots,x_n)| x_i\geq 0,  \sum x_i =1\}$ be the standard simplex. 
Consider the subset  $C\subset S_n$ consisting of functions with exactly $n+1$ intervals of constant sign. The lengths of these intervals are positive numbers $x_0,x_1,\dots,x_n$ satisfying $\sum x_i =1$, and a function from $C$ is determined by $x=(x_0,\dots,x_n)$ and the sign $\pm$ that the function has on the first interval. Thus we obtain two embeddings $\psi^n_{\pm}: {\rm Int}\ \Delta^n\to C$, and  $C$ is the disjoint union of the images of $\psi^n_+$ and $\psi^n_-$.

The maps $\psi^n_{\pm}$ extend continuously to the boundary $\partial \Delta^n$: when some $x_i$s shrink to zero, the respective segments of constant sign of a function disappear, and if the function has the same signs in the neighboring segments, they merge together. For example, let $n=2$. Then $\psi^2_+ (0,x_1,x_2)$
has two intervals of constant sign and equals $\psi^1_- (x_1,x_2)$, whereas $\psi^2_+ (x_0,0,x_2)$ is constant function with value $+1$, i.e., equals $\psi^0_+(1)$.

We have: $S_n - C = S_{n-1}$, and the latter is homeomorphic to $S^{n-1}$ by the induction assumption.  Each  map $\psi^n_{\pm}$ sends $\partial \Delta^n$  to $S_{n-1}$, and we claim that the degree of $\psi^n_{\pm}$ is one. Indeed, the faces of $\partial \Delta^n$ are given by one of the conditions: $x_0=0, x_1=0,\dots,x_n=0$. Since 
$\psi^n_{\pm}(0,x_1,\dots,x_n)=\psi^{n-1}_{\mp}(x_1,\dots,x_n)$ and    
$\psi^n_{\pm}(x_0,\dots,x_{n-1},0)=\psi^{n-1}_{\pm}(x_0,\dots,x_{n-1}),$ 
the map $\psi^n_{\pm}$ sends the faces $x_0=0$ and $x_n=0$ to the two $n-1$-dimensional cells of $S_{n-1}$, and the other faces are sent  to the $n-2$-skeleton of $S_{n-1}$. Therefore deg $\psi^n_{\pm}=1$.

Since the attaching maps of two $n$-dimensional discs $\Delta^n$ to $S^{n-1}$ have degree one,  $S_n$ is $n$-dimensional sphere.
\proofend

One can also consider  the space  of piece-wise constant function on the circle with values $\pm 1$ and at most $n$ intervals of constant sign ($n$ even). Such a space is also homeomorphic to $S^{n}$: cut the circle at, say, point $0$ to obtain a piece-wise constant function on an interval with at most $n+1$
intervals of constant sign, and apply Theorem \ref{sphere}.

\section{Proof of the converse Ghys theorem} \label{Ghys}

Let us start with a reformulation described in \cite{O-T1}.

A diffeomorphism $f:\RP^1\to\RP^1$ has a unique lifting to a homogeneous of degree one area preserving diffeomorphism $F$ of the punctured plane. If $f$ is a projective transformation then $F \in SL(2,\R)$. Let $x$ be the angular parameter on $\RP^1$ so that $x$ and $x+\pi$ describe the same point. Then $(x,r)$ are the polar coordinates in the plane and
$$
F(x,r)=(f(x),rf'^{-1/2}(x)).
$$
Let $\gamma(x)$ be the image of the unit circle under $F$, this is a centrally symmetric curve that bounds area $\pi$. The curve $\gamma$ satisfies the differential equation
\begin{equation} \label{diffeq}
\gamma''(x)=-k(x)\gamma(x)
\end{equation}
 where $k(x)$ is a $\pi$-periodic function called the potential. The relation of the potential with the Schwarzian derivative is as follows:
$$
k=\frac{1}{2}S(f)+1.
$$
In particular, the zeros of the Schwarzian corresponds to the values 1 of the function $k(x)$ (indeed, if $k(x)\equiv 1$ then $\gamma$ is a central ellipse, $F\in SL(2,\R)$ and $f$ is a projective transformation).

Thus we arrive at the following reformulation of Theorem \ref{Schw}: {\it if a function $k(x)-1$ on $\RP^1$ changes sign at least four times then there exists an orientation preserving diffeomorphism $\varphi$ of the projective line such that the function ${\bar k}=\varphi^*(k)$ is the potential of a centrally symmetric closed parametric curve $\gamma(x)$ in the punctured plane bounding area $\pi$, that is, a curve satisfying the differential equation $\gamma''(x)=-{\bar k}(x)\gamma(x)$.}

The proof consists of the same four steps as in Example \ref{exa}. 

\paragraph{\it Step 1.}
Let $k_1,k_2$ be two positive numbers satisfying $k_1>1, k_1+k_2=2$ and both sufficiently close to $1$. We claim that there exists a $\pi$-periodic step function $h(x)$  with four intervals of constant values $k_1,k_2,k_1,k_2$ on $[0,\pi]$  such that the respective solution of the differential equation (\ref{diffeq}) is a closed curve.

To prove this, consider the frame $F(x)=(\gamma(x),\gamma'(x))$. The differential equation (\ref{diffeq}) rewrites as
\begin{equation} \label{frame}
F'(x)=  F(x) A(x)
\end{equation}
where 
$$
A(x)=\begin{pmatrix}
0 & -k(x)\\
1 & 0
\end{pmatrix}.
$$
Equation (\ref{frame}) defines a curve on the group $SL(2,\R)$; the curve $\gamma$ is centrally symmetric and closed iff $F(\pi)=-F(0)$. Let us refer to the last equality as the monodromy condition.

Let the desired step function $h(x)$ have intervals of constant values of lengths $t_1,t_2,t_3,t_4$ with $t_1+t_2+t_3+t_4=\pi$. For a constant potential $k$, equation (\ref{frame}) is easily solved:
$$
F(x)=F(0)e^{xA}=F(0)\begin{pmatrix}
\cos(\sqrt{k}x) & -\sqrt{k}\sin(\sqrt{k}x)\\
\frac{1}{\sqrt{k}}\sin(\sqrt{k}x) & \cos(\sqrt{k}x)
\end{pmatrix}.
$$
It follows that the monodromy condition is
\begin{equation} \label{mono}
e^{t_1A} e^{t_2B} e^{t_3A} e^{t_4B} = -E
\end{equation}
where
$$
A=\begin{pmatrix}
0 & -k_1\\
1 & 0
\end{pmatrix}, \ \ 
B=\begin{pmatrix}
0 & -k_2\\
1 & 0
\end{pmatrix}
$$
and $E$ is the unit matrix.

Let us look for a solution satisfying $t_3=t_1, t_4=t_2$; then $t_1+t_2=\pi/2$. Set: $\alpha=t_1\sqrt{k_1}, \beta=t_2\sqrt{k_2}$. A direct computation shows that (\ref{mono}) is satisfied once 
\begin{equation} \label{tan}
\tan \alpha \tan \beta = \sqrt{k_1k_2}.
\end{equation}
 The constraint on $\alpha$ and $\beta$ is
$$
\frac{\alpha}{\sqrt{k_1}}+\frac{\beta}{\sqrt{k_2}}=\frac{\pi}{2}.
$$
If $\alpha$ is close to $\pi/2$ then the left hand side of (\ref{tan}) is greater, and if $\alpha$ is close to $0$ then it is smaller than the right hand side. It follows that (\ref{tan}) has a solution.

\paragraph{\it Step 2.}
Since $k(x)-1$ changes sign at least four times, there is a constant $c>0$ such that $k$ takes the values $1+c,1-c,1+c,1-c$ at points, say, $x_1,x_2,x_3,x_4$. Let $k_1=1+c, k_2=1-c$, and let $h(x)$ be the step function from Step 1. For every $\eps >0$, there exists a diffeomorphism $\varphi \in {\rm Diff}_+ (\RP^1)$ which stretches neighborhoods of the points  $x_1,\dots,x_4$ so that the function $\varphi^*(k)$ is $\eps$-close in measure to $h$.

\paragraph{\it Step 3.}
Similarly to Step 3 in Section \ref{SHK}, consider a 3-parameter family of diffeomorphisms $\psi_{\alpha}\in {\rm Diff} (\RP^1)$ that change the intervals of constant values of the step function $h(x)$. Given $\alpha$, consider the function $\psi_{\alpha}^*(h)$ as the potential of equation (\ref{frame}) with the initial conditions $F(0)=E$. The formula $G(\alpha)=F(\pi)$ defines a smooth map $D^3 \to SL(2,\R)$ that takes the origin to the matrix $-E$.

\begin{lemma} \label{nondeg}
The differential $dG$ is non-degenerate at the origin.
\end{lemma} 

\paragraph{\bf Proof} 
Stretch the intervals of constant values of the potential function to $t_i+\eps s_i,\ i=1,2,3,4$; the vector $s=(s_1,s_2,s_3,s_4),\ s_1+s_2+s_3+s_4=0$ is interpreted as a tangent vector to $D^3$ at the origin. Using the formula for monodromy (\ref{mono}), we compute:
\begin{equation} \label{onto}
-dG(s)=s_1 A+s_4 B + s_2 e^{t_1A}Be^{-t_1A}+s_3 e^{t_2B}Ae^{t_2B}
\end{equation}
where $A,B,t_1,t_2$ are as in Step 1. We need to check that the linear map $dG:\R^4\to sl_2$, given by (\ref{onto}), is surjective and that its kernel is transverse to the hyperplane $s_1+s_2+s_3+s_4=0$. Both claims follow, by a direct computation, from the explicit formulas for the matrices $A,B$ and their exponents given in Step 1.
\proofend

\paragraph{\it Step 4.}
This last step is identical to Step 4 in Section \ref{SHK}: replace the potential $h$ in the definition of the map $G$ in Step 3 by $\varphi^*(k)$. We obtain a new monodromy map ${\bar G}: D^3\to SL(2,\R)$ whose image contains the matrix $-E$. The respective curve closes up, and we are done.

\begin{remark} \label{hyp}
{\rm The Ghys theorem is closely related to the four vertex theorem in the hyperbolic plane \cite{Sin}. Let $\gamma$ be an oval in $H^2$. Each tangent line to $\gamma$ intersects the circle at infinity  at two points, and this defines a circle diffeomorphism $f_{\gamma}$. In the projective model of hyperbolic geometry, the circle at infinity is represented by a conic in $\RP^2$. A conic has a canonical projective structure, hence $f_{\gamma}$ can be viewed as a diffeomorphism of $\RP^1$. Singer's theorem asserts that the zeros of the Schwarzian $S(f_{\gamma})$ correspond to the vertices of $\gamma$ (in the hyperbolic metric, of course), see \cite{O-T} for a discussion. 

Note however that  a converse four vertex theorem for the hyperbolic plane does not hold in the same way  as in the Euclidean plane: if the positive curvature function is too small then the respective curve in the hyperbolic plane does not close up.
}
\end{remark}

\section{Problems and conjectures} \label{conj}

There are many other results extending the four vertex theorem. In each case, it is interesting to find the  converse theorem; we mention but a few.

\paragraph{\it Problem 1.} Another classic theorem of Mukhopadhyaya is that a plane oval has at least six affine vertices (also known as sextactic points). An affine vertex is a point at which the curve is abnormally well approximated by a conic: at a generic point, a conic passes through five infinitesimally close points of the curve, whereas at an affine point, this number equals six. Every oval $\gamma$ can be given an {\it affine parameterization} such that $\det (\gamma'(x),\gamma''(x))$ is constant. Then 
$\gamma'''(x)=-k(x)\gamma'(x)$ where the function $k(x)$ is called the affine curvature. The affine vertices are the critical points of the affine curvature, see, e.g., \cite{O-T}.

A conjectural converse theorem asserts that {\it if a periodic function $k(x)$ has at least six extrema then there exists a plane oval $\gamma(x)$ whose affine curvature at point $\gamma(x)$ is $k(x)$} (of course, here $x$ is not necessarily an affine parameter).

\paragraph{\it Problem 2.} The four vertex theorem has numerous discrete versions, see, e.g.,  \cite{O-T,Pak} for surveys and references. For example, let $P$ be a convex $n$-gon with vertices $x_1,\dots,x_n$. Assume that $n\geq 4$ and that no four consecutive vertices lie on a circle. Consider the circles circumscribing  triples of  consecutive vertices $x_{i-1} x_i x_{i+1}$, and assume that the center of this circle lies inside the cone of the vertex $x_i$ (such a polygon is called {\it coherent}). Let $r_1,\dots,r_n$ be the cyclic sequence of the radii of the circles. Then the sequence $r_1,\dots,r_n$ has at least two local maxima and two local minima.  

A conjectural converse theorem asserts that {\it if a  cyclic sequence $r_1,\dots,r_n$ has at least two local maxima and two local minima then it corresponds, as described above, to a coherent convex polygon.}

Another version of discrete four vertex theorem concerns the circles tangent to the triples of consecutive sides of a polygon: the radii of such inscribed circles also form a cyclic sequence with at least two local maxima and two local minima. One conjectures that a converse theorem holds as well. 

\bigskip

{\bf Acknowledgments}. It is a pleasure to acknowledge interesting discussions with D. Fuchs, H. Gluck, D. Khavinson, V. Ovsienko, I. Pak, A. Pinkus, R. Schwartz, D. Singer and  V. Totik.  The author  was partially supported by an NSF grant DMS-0555803.

\end{document}